\documentclass[11pt]{amsart}

\def\NZQ{\mathbb}               

\def\QQ{{\NZQ Q}}
\def\ZZ{{\NZQ Z}}
\def\RR{{\NZQ R}}

\def\frk{\mathfrak}               

\def\Phi{{\frk N}}
\def\a{\alpha}
\def\b{\beta}

\def\d{\delta}
\def\opn#1#2{\def#1{\operatorname{#2}}} 
\opn\chara{char} \opn\length{\ell} \opn\pd{pd} \opn\rk{rk}
\opn\projdim{proj\,dim} \opn\injdim{inj\,dim} \opn\rank{rank}
\opn\depth{depth} \opn\grade{grade} \opn\height{height}
\opn\embdim{emb\,dim} \opn\codim{codim}

\opn\Tr{Tr} \opn\bigrank{big\,rank}
\opn\superheight{superheight}\opn\lcm{lcm}
\opn\trdeg{tr\,deg}
\opn\reg{reg} \opn\lreg{lreg} \opn\ini{in} \opn\lpd{lpd}
\opn\size{size}\opn{\mult}{mult}
\opn\div{div} \opn\Div{Div} \opn\cl{cl} \opn\Cl{Cl}
\opn\Spec{Spec} \opn\Supp{Supp} \opn\supp{supp} \opn\Sing{Sing}
\opn\Ass{Ass} \opn\Min{Min}
\opn\Ann{Ann} \opn\Rad{Rad} \opn\Soc{Soc}
\opn\Syz{Syz} \opn\Im{Im} \opn\Ker{Ker} \opn\Coker{Coker}
\opn\Am{Am} \opn\Hom{Hom} \opn\Tor{Tor} \opn\Ext{Ext}
\opn\End{End} \opn\Aut{Aut} \opn\id{id}

\opn\nat{nat}
\opn\pff{pf}
\opn\Pf{Pf} \opn\GL{GL} \opn\SL{SL} \opn\mod{mod} \opn\ord{ord}
\opn\Gin{Gin}
\opn\Hilb{Hilb}\opn\adeg{adeg}\opn\std{std}\opn\ip{infpt}
\opn\Pol{Pol}
\opn\sat{sat}
\opn\Var{Var}
\opn\Gen{Gen}
\opn\vol{vol}

\opn\aff{aff} \opn\con{conv} \opn\relint{relint} \opn\st{st}
\opn\lk{lk} \opn\cn{cn} \opn\core{core} \opn\vol{vol}
\opn\link{link} \opn\star{star}
\opn\gr{gr}

\def\Cc{{\mathcal C}}

\def\Fc{{\mathcal F}}
\def\Pc{{\mathcal P}}

\def\pot#1#2{#1[\kern-0.28ex[#2]\kern-0.28ex]}

\opn\dirlim{\underrightarrow{\lim}}
\opn\inivlim{\underleftarrow{\lim}}
\def\Implies{\ifmmode\Longrightarrow \else
        \unskip${}\Longrightarrow{}$\ignorespaces\fi}
\def\implies{\ifmmode\Rightarrow \else
        \unskip${}\Rightarrow{}$\ignorespaces\fi}
\def\iff{\ifmmode\Longleftrightarrow \else
        \unskip${}\Longleftrightarrow{}$\ignorespaces\fi}

\let\:=\colon
\newtheorem{Theorem}{Theorem}[section]
\newtheorem{Lemma}[Theorem]{Lemma}

\newtheorem{Example}[Theorem]{Example}

\let\epsilon\varepsilon
\let\phi=\varphi
\let\kappa=\varkappa
\def\qed{\ifhmode\textqed\fi
      \ifmmode\ifinner\quad\qedsymbol\else\dispqed\fi\fi}
\def\textqed{\unskip\nobreak\penalty50
       \hskip2em\hbox{}\nobreak\hfil\qedsymbol
       \parfillskip=0pt \finalhyphendemerits=0}
\def\dispqed{\rlap{\qquad\qedsymbol}}

\opn\dis{dis}
\def\pnt{{\raise0.5mm\hbox{\large\bf.}}}

\opn\Lex{Lex}

\begin{document}

\title{Ehrhart polynomials of convex polytopes
with small volumes}

\author{Takayuki Hibi, Akihiro Higashitani
and Yuuki Nagazawa}

\address{Takayuki Hibi,
Department of Pure and Applied Mathematics,
Graduate School of Information Science and Technology,pe
Osaka University,
Toyonaka, Osaka 560-0043, Japan}
\email{hibi@math.sci.osaka-u.ac.jp}

\address{Akihiro Higashitani,
Department of Pure and Applied Mathematics,
Graduate School of Information Science and Technology,
Osaka University,
Toyonaka, Osaka 560-0043, Japan}
\email{sm5037ha@ecs.cmc.osaka-u.ac.jp}

\address{Yuuki Nagazawa,
Department of Pure and Applied Mathematics,
Graduate School of Information Science and Technology,
Osaka University,
Toyonaka, Osaka 560-0043, Japan}
\email{sm5032ny@ecs.cmc.osaka-u.ac.jp}

\begin{abstract}
Let $\Pc \subset \RR^d$ be an integral
convex polytope of dimension $d$ and
$\delta(\Pc)
= (\delta_0, \delta_1, \ldots, \delta_d)$
be its $\delta$-vector.
By using the known inequalities on
$\delta$-vectors,
we classify the possible $\delta$-vectors
of convex polytopes of dimension $d$
with $\sum_{i=0}^{d} \delta_i \leq 3$.
\end{abstract}
\maketitle

\section*{Introduction}
One of the most attractive problems on enumerative
combinatorics of convex polytopes is to find
a combinatorial characterizations of the Ehrhart
polynomials of integral convex polytopes.
First of all, we recall what the Ehrhart polynomial
of a convex polytope is.

Let $\Pc \subset \RR^N$ be an {\em integral}
convex polytope; i.e., a convex polytope
any of whose vertices has integer coordinates,
of dimension $d$, and let $\partial \Pc$ denote
the boundary of $\Pc$.
Given a positive integer $n$ we define
the numerical functions $i(\Pc,n)$ and
$i^*(\Pc,n)$ by setting
\[
i(\Pc,n) = |n\Pc \cap \ZZ^N|, \, \, \, \, \,
i^*(\Pc,n) = |n(\Pc - \partial \Pc) \cap \ZZ^N|.
\]
Here $n\Pc = \{ n\alpha : \alpha \in \Pc \}$
and $|X|$ is the cardinality of a finite set $X$.

The systematic study of $i(\Pc,n)$ originated
in the work of Ehrhart \cite{Ehrhart}, who established
the following fundamental properties:
\begin{enumerate}
\item[(0.1)]
$i(\Pc,n)$ is a polynomial in $n$ of degree $d$;
(Thus in particular $i(\Pc,n)$ can be defined
for {\em every} integer $n$.)
\item[(0.2)]
$i(\Pc,0) = 1$;
\item[(0.3)]
(loi de r\'eciprocit\'e)
$i^*(\Pc,n) = ( - 1 )^d i(\Pc, - n)$
for every integer $n > 0$.
\end{enumerate}
We say that $i(\Pc,n)$ is the {\em Ehrhart polynomial}
of $\Pc$.
An introduction
to the theory of Ehrhart polynomials
is discussed in
\cite[pp. 235--241]{StanleyEC}
and \cite[Part II]{HibiRedBook}.

We define the sequence
$\delta_0, \delta_1, \delta_2, \ldots$ of integers
by the formula
\begin{eqnarray}
\label{delta}
(1 - \lambda)^{d + 1}
\left[ 1 + \sum_{n=1}^{\infty} i(\Pc,n) \lambda^n \right]
= \sum_{i=0}^{\infty} \delta_i \lambda^i.
\end{eqnarray}
Then the basic fact (0.1) and (0.2) on $i(\Pc,n)$
together with a fundamental result on generating
function (\cite[Corollary 4.3.1]{StanleyEC})
guarantee that $\delta_i = 0$ for every $i > d$.
We say that the sequence
\[
\delta(\Pc) = (\delta_0, \delta_1, \ldots, \delta_d)
\]
which appears in Eq.\,(\ref{delta}) is the
{\em $\delta$-vector} of $\Pc$.
Thus $\delta_0 = 1$ and
$\delta_1 = |\Pc \cap \ZZ^N| - (d + 1)$.

It follows from loi de r\'eciprocit\'e (0.3) that
\begin{eqnarray}
\label{deltadual}
\sum_{n=1}^{\infty} i^*(\Pc,n) \lambda^n
= \frac
{\sum_{i=0}^{d} \delta_{d-i} \lambda^{i+1}}
{(1 - \lambda)^{d + 1}}.
\end{eqnarray}
In particular
$\delta_d = |(\Pc - \partial \Pc) \cap \ZZ^N|$.
Hence $\delta_1 \geq \delta_d$.
Moreover, each $\delta_i$ is nonnegative
(Stanley \cite{StanleyDRCP}).
In addition,
if $(\Pc - \partial \Pc) \cap \ZZ^N$ is nonempty,
i.e., $\delta_d \neq 0$, then one has
$\delta_1 \leq \delta_i$ for every $1 \leq i < d$
(\cite{HibiLBT}).
When $d = N$, the leading coefficient
$(\sum_{i=0}^{d}\delta_i)/d!$
of $i(\Pc,n)$ is equal to the usual volume
$\vol(\Pc)$ of $\Pc$
(\cite[Proposition 4.6.30]{StanleyEC}).

It follows from Eq.\,(\ref{deltadual}) that
\[
\max\{ i : \delta_i \neq 0 \}
+ \min\{ i :
i(\Pc - \partial \Pc) \cap \ZZ^N \neq \emptyset \}
= d + 1.
\]
Let $s = \max\{ i : \delta_i \neq 0 \}$.
Stanley \cite{StanleyJPAA} shows the inequalities
\begin{eqnarray}
\label{Stanley}
\delta_0 + \delta_1 + \cdots + \delta_i
\leq \delta_s + \delta_{s-1} + \cdots + \delta_{s-i},
\, \, \, \, \,
0 \leq i \leq [s/2]
\end{eqnarray}
by using the theory of Cohen--Macaulay rings.
On the other hand, the inequalities
\begin{eqnarray}
\label{Hibi}
\delta_{d-1} + \delta_{d-2} + \cdots + \delta_{d-i}
\leq \delta_2 + \delta_3 + \cdots + \delta_i
+ \delta_{i+1},
\, \, \, \, \,
1 \leq i \leq [(d-1)/2]
\end{eqnarray}
appear in \cite[Remark (1.4)]{HibiLBT}.

Somewhat surprisingly, when $\sum_{i=0}^{d} \delta_i \leq 3$,
the above inequalities
(\ref{Stanley}) together with (\ref{Hibi})
give a characterization
of the possible $\delta$-vectors.  In fact,

\begin{Theorem}
\label{main}
Let $d \geq 3$.
Given a finite sequence
$(\delta_0, \delta_1, \ldots, \delta_d)$
of nonnegative integers, where $\delta_0 = 1$
and $\delta_1 \geq \delta_d$,
which satisfies $\sum_{i=0}^{d} \delta_i \leq 3$,
there exists an integral convex
polytope $\Pc \subset \RR^d$
of dimension $d$
whose $\delta$-vector
coincides with
$(\delta_0, \delta_1, \ldots, \delta_d)$
if and only if
$(\delta_0, \delta_1, \ldots, \delta_d)$
satisfies all inequalities
{\em (\ref{Stanley})} and {\em (\ref{Hibi})}.
\end{Theorem}

The ``Only if'' part of Theorem \ref{main}
is obvious.  In addition, no discussion
will be required for the case of
$\sum_{i=0}^{d} \delta_i = 1$.
The ``If'' part of Theorem \ref{main} will be given
in Section $2$ for the case of
$\sum_{i=0}^{d} \delta_i = 2$
and in Section $3$ for the case of
$\sum_{i=0}^{d} \delta_i = 3$.

On the other hand, 
Example \ref{higashitani}
shows that
Theorem \ref{main} is no longer true
for the case of $\sum_{i=0}^{d} \delta_i = 4$.
Finally,
when $d \leq 2$,
the possible $\delta$-vectors
are known (Scott \cite{Scott}).

\section{Review on the computation of the $\delta$-vector of a simplex}
We recall from \cite[Part II]{HibiRedBook}
the well-known combinatorial technique how to compute
the $\delta$-vector of a simplex.

\begin{itemize}
\item
Given an integral $d$-simplex $\Fc \subset \RR ^N$ with the vertices
$v_0, v_1, \ldots, v_d$, we set
$\widetilde \Fc=\left\{(\a,n)\in \RR^{N+1} \, : \, \a \in \Fc \right\}$,
which is an integral $d$-simplex in $\RR^{N+1}$
with the vertices
$(v_0,1), (v_1,1), \ldots,(v_d,1)$
and
$
\partial \widetilde \Fc=\left\{(\a,1)\in \RR^{N+1} \, : \, \a \in \partial \Fc \right\}
$
is its boundary.
Clearly
$
i(\Fc,n)=i(\widetilde \Fc,n)
$
and
$
i^*(\Fc,n)=i^*(\widetilde \Fc,n)
$
for all $n$.
\item
The subset $\Cc=\Cc(\widetilde \Fc) \subset \RR^{N+1}$ defined by
$
\Cc=\{r\b \, : \, \b \in \widetilde \Fc,0 \leq r \in \QQ \}
$
is called \textit{the simplicial cone associated with $\Fc \subset \RR^N$ with apex $(0, \ldots ,0)$}.
Its boundary is
$
\partial \Cc = \left\{r\b \, : \,  \b \in \partial \widetilde \Fc,0 \leq r \in \QQ \right\}.
$
One has
$
i(\Fc,n)= | \left\{ (\a,n) \in \Cc \, : \,  \a \in \ZZ^N \right\} |
$ and
$
i^*(\Fc,n)= | \left\{ (\a,n) \in \Cc \setminus \partial \Cc \, : \,  \a \in \ZZ^N \right\} |.
$
\item
Each rational point $\a \in \Cc$
has a unique expression of the form
$
\a= \sum_{i=0}^{d}r_i(v_i,1)
$
with each $0 \leq r_i \in \QQ$.
Moreover,
each rational point $\a \in \Cc \setminus \partial \Cc$
has a unique expression of the form
$
\a= \sum_{i=0}^{d}r_i(v_i,1)
$
with each $0<r_i \in \QQ$.
\item
Let $S$ (resp. $S^*$) be the set of all points $\a \in \Cc \cap \ZZ^{N+1}$
(resp. $\a \in (\Cc- \partial \Cc) \cap \ZZ^{N+1}$)
of the form
$
\a= \sum_{i=0}^{d}r_i(v_i,1),
$
where each $r_i \in \QQ$ with $0 \leq r_i<1$ (resp. with $0<r_i \leq 1$).
\item
The degree of an integer point $(\a,n) \in \Cc$ is
$
\deg(\a,n):=n.
$
\end{itemize}

\begin{Lemma}
\label{computation}
{\em (a)}
Let $\delta_i$ be the number of integer points $\a \in S$
with $\deg \a=i$.  Then
\[
1+ \sum_{n=1}^{\infty}i(\Fc,n) \lambda^n=\frac{\delta_0+\cdots+\delta_d\lambda^d}{(1-\lambda)^{d+1}}.
\]

{\em (b)}
Let $\delta_i^*$ be the number of integer points $\a \in S^*$
with $\deg\a=i$.  Then
\[
\sum_{n=1}^{\infty}i^*(\Fc,n) \lambda^n=\frac{\delta_1^*\lambda+\cdots+\delta_{d+1}^*\lambda^{d+1}}{(1-\lambda)^{d+1}}.
\]

{\em (c)}
One has $\delta_i^* = \delta_{(d+1)-i}$ for each $1 \leq i \leq d+1$.
\end{Lemma}

\begin{Example}
\label{higashitani}
{\em
Theorem \ref{main} is no longer true
for the case of $\sum_{i=0}^{d} \delta_i = 4$.
In fact, 
the sequence $(1,0,1,0,1,1,0,0)$
cannot be the $\delta$-vector of 
an integral convex polytope of dimension $7$.
Suppose, on the contrary, 
that there exists an integral convex polytope 
$\Pc \subset \RR^N$ with 
$(\delta_0, \delta_1, \ldots, \delta_7)=(1,0,1,0,1,1,0,0)$ 
its $\delta$-vector.
Since $\delta_1 = 0$,
we know that
$\Pc$ is a simplex.
Let $v_0,v_1,\ldots,v_7$ be the vertices of $\Pc$.
By using Lemma \ref{computation}, one has
$S = \{(0,\ldots,0),(\alpha,2),(\beta,4),(\gamma,5)\}$ 
and
$S^* = \{(\alpha',3),(\beta',4),(\gamma',6),
(\sum_{i=0}^7v_i,8)\}$.
Write $\alpha'=\sum_{i=0}^{7}r_iv_i$
with each $0 < r_i \leq 1$.
Since $(\alpha',3) \not\in S$,
there is $0 \leq j \leq 7$ with $r_j = 1$. 
If there are $0 \leq k < \ell \leq 7$
with $r_k = r_\ell = 1$,
say, $r_0 = r_1 = 1$, 
then $0 < r_q < 1$ for each 
$2 \leq q \leq 7$ and 
$\sum_{i=2}^{7} r_i = 1$.
Hence $(\alpha' - v_0 - v_1, 1) \in S$,
a contradiction.
Thus there is a unique 
$0 \leq j \leq 7$ with $r_j = 1$,
say, $r_0 = 1$. 
Then $\alpha=\sum_{i=1}^{7}r_iv_i$ and
$\gamma=\sum_{i=1}^{7}(1-r_i)v_i$.
Let $\Fc$ denote the facet of $\Pc$ whose vertices are $v_1,v_2,\ldots,v_7$
with $\delta(\Fc) = (\delta'_0, \delta'_1, \ldots, \delta'_6) \in \ZZ^7$.
Then $\delta'_2 = \delta'_5 = 1$.  
Since $\delta'_i \leq \delta_i$ for each $0 \leq i \leq 6$,
it follows that $\delta(\Fc) = (1,0,1,0,0,1,0)$. 
This contradicts the inequalities (\ref{Stanley}). 
}
\end{Example}

\section{A proof of Theorem \ref{main} when
$\sum_{i=0}^{d} \delta_i = 2$}
The goal of this section is to prove the ``If'' part
of Theorem \ref{main} when $\sum_{i=0}^{d} \delta_i = 2$.
First of all, we recall the following well-known

\begin{Lemma}
\label{wellknown}
Suppose that $(\delta_0, \delta_1, \ldots, \delta_d)$ is the $\delta$-vector
of an integral convex polytope of dimension $d$.
Then there exists
an integral convex polytope of dimension $d + 1$ whose $\delta$-vector
is $(\delta_0, \delta_1, \ldots, \delta_d, 0)$.
\end{Lemma}

Let $d \geq 3$.
We study a finite sequence
$(\delta_0, \delta_1, \ldots, \delta_d)$
of nonnegative integers with $\delta_0 = 1$
and $\delta_1 \geq \delta_d$
which satisfies all inequalities
(\ref{Hibi})
together with
$\sum_{i=0}^{d} \delta_i =2$.
Since $\delta_0 = 1$,
$\delta_1 \geq \delta_d$
and $\sum_{i=0}^{d} \delta_i =2$,
one has $\delta_d = 0$.
Hence
there is an integer
$i \in \{ 1, \ldots, \left[(d-1)/ 2 \right] +1 \}$
such that
$(\delta_0, \delta_1, \ldots, \delta_d)
= (1,0,\ldots,0,\underbrace{1}_{i\text{-th}},0,\ldots,0)$,
where $\underbrace{1}_{i\text{-th}}$ stands for $\delta_i=1$.
By virtue of Lemma \ref{wellknown}
our work is to find an integral convex polytopes $\Pc$
of dimension $d$
with $(1,0,\ldots,0,\underbrace{1}_{((d-1)/2)\text{-th}},0,\ldots,0)
\in \ZZ^{d+1}$
its $\delta$-vector.

Let $\Pc \subset \RR^d$ be the integral simplex of dimension $d$
whose vertices $v_0, v_1, \ldots, v_d$ are
\[
v_i =
\begin{cases}
(0,\ldots,0,\underbrace{1}_{i\text{-th}},\underbrace{1}_{(i+1)\text{-th}},0,\ldots,0), &i=1,\ldots,d-1, \\
(1,0,\ldots,0,1), &i=d, \\
(0,0,\ldots\ldots,0), &i=0.
\end{cases}
\]
When $d$ is odd, one has $\vol(\Pc) = 2/d!$ by using an elementary linear algebra.
Since
\[
\frac{1}{2}\left\{ (v_0,1)+(v_1,1)+\cdots+(v_d,1)\right\} = (1,1,\ldots,1,(d-1)/2) \in \ZZ^{d+1},
\]
Lemma\ref{computation} says that
$\delta_{(d-1)/1} \geq 1$.
Thus, since $\vol(\Pc)=2/d!$,
one has
\[\delta(\Pc) =
(1,0,\ldots,0,\underbrace{1}_{((d-1)/2)\text{-th}},0,\ldots,0),
\]
as desired.

\section{A proof of Theorem \ref{main} when
$\sum_{i=0}^{d} \delta_i = 3$}
The goal of this section is to prove the ``If'' part
of Theorem \ref{main} when $\sum_{i=0}^{d} \delta_i = 3$.
Let $d \geq 3$.
Suppose that a finite sequence
$(\delta_0, \delta_1, \ldots, \delta_d)$
of nonnegative integers with $\delta_0 = 1$
and $\delta_1 \geq \delta_d$
satisfies all inequalities
(\ref{Stanley}) and (\ref{Hibi})
together with
$\sum_{i=0}^{d} \delta_i =3$.

When there is $1 \leq i \leq d$ with $\delta_i = 2$,
the same discussion as in Section $1$ can be applied.
In fact, instead of the vertices of the convex polytope
arising in the last paragraph of Section $1$,
we may consider the convex polytope
whose vertices $v_0, v_1, \ldots, v_d$ are
\begin{eqnarray*}
v_i:=
\begin{cases}
(0,\ldots,0,\underbrace{1}_{i\text{-th}},\underbrace{1}_{(i+1)\text{-th}},0,\ldots,0), &i=1,\ldots,d-1, \\
(2,0,\ldots,0,1), &i=d, \\
(0,0,\ldots\ldots,0), &i=0.
\end{cases}
\end{eqnarray*}

Now, in what follows, a sequence
$(\delta_0, \delta_1, \ldots, \delta_d)$
with each $\delta_i \in \{ 0, 1 \}$,
where
$\delta_0 = 1$ and $\delta_1 \geq \delta_d$,
which satisfies all inequalities
(\ref{Stanley}) and (\ref{Hibi})
together with
$\sum_{i=0}^{d} \delta_i = 3$
will be considered.

If $\delta_d = 1$, then $\delta_1 = 1$.  However,
since $d \geq 3$, this contradicts (\ref{Stanley}).
If $\delta_1 = 1$, then $\delta_2 = 1$ by (\ref{Stanley}).
Clearly,
$(1,1,1,0,\ldots,0) \in \ZZ^{d+1}$
is a possible $\delta$-vectors.
Thus we will assume that $\delta_1 = \delta_d = 0$.
Let $\d_m=\d_n=1$ with  $1< m < n < d$.
Let $p=m-1$, $q=n-m-1$, and $r=d-n$.
By (\ref{Stanley}) one has $0 \leq q \leq p$ .
Moreover, by (\ref{Hibi}) one has $p \leq r$.
Consequently,
\begin{eqnarray}
\label{hhn}
0 \leq q \leq p \leq r, \, \, \, \, \, p+q+r=d-2.
\end{eqnarray}
Our work is to construct an integral convex polytope $\Pc$ with dimension $d$
whose $\d$-vector coincides with
$\d(\Pc)=(1,\underbrace{0,\ldots,0}_{p},1,\underbrace{0,\ldots,0}_{q},1,\underbrace{0,\ldots,0}_{r})$
for an arbitrary integer $1< m < n < d$ satisfying the conditions (\ref{hhn}).


\begin{Lemma}
\label{first}
Let $d=3k+2$.  There exists an integral convex polytope $\Pc$ of dimension $d$
whose $\d$-vector coincides with
\[
(1,\underbrace{0,\ldots,0}_{k},1,\underbrace{0,\ldots,0}_{k},1,\underbrace{0,\ldots,0}_{k})
\in \ZZ^{d+1}.
\]
\end{Lemma}

\begin{proof}
When $k \geq 1$, let $\Pc \subset \RR^d$ be the integral simplex of dimension $d$
with the vertices
$v_0, v_1, \ldots, v_d$, where
\[
v_i =
\begin{cases}
(0,\ldots,0,
\underbrace{1}_{i\text{-th}},
\underbrace{1}_{(i+1)\text{-th}},
\underbrace{1}_{(i+2)\text{-th}},0,\ldots,0), \qquad &i=1,\ldots,d-2, \\
(1,0,\ldots,0,1,1), &i=d-1, \\
(1,1,0,\ldots,0,1), &i=d, \\
(0,\ldots \ldots \ldots ,0), &i=0.
\end{cases}
\]
By using the induction on $k$
it follows that $\vol(\Pc) = 3 / d!$.
Since
\begin{eqnarray*}
\frac{1}{3}\left\{ (v_0,1)+(v_1,1)+\cdots+(v_d,1)\right\}
= (1,1,\ldots,1,k+1) \in \ZZ^{d+1},
\end{eqnarray*}
Lemma \ref{computation} now guarantees that
$\d_{k+1}\geq 1$ and $\d_{k+1}^* \geq 1$.
Hence $\d_{k+1} = 1$ and $\d_{2k+2}=1$, as required.
\hspace{12.15cm}
\end{proof}


\begin{Lemma}
\label{second}
Let $d=3k+2$, $\ell >0$ and $d'=d+2\ell$.
There exists an integral simplex $\Pc \subset \RR^{d'}$
of dimension $d'$ whose $\d$-vector coincides with
\[
(1,\underbrace{0,\ldots,0}_{k+\ell},1,\underbrace{0,\ldots,0}_{k},1,\underbrace{0,\ldots,0}_{k+\ell})
\in \ZZ^{d'+1}.
\]
\end{Lemma}

\begin{proof}
{\bf (First Step)}
Let $k=0$.  Thus $d=2$ and $d'=2\ell+2$.
Let $\Pc \subset \RR^{d'}$ be an integer convex polytope 
of dimension $d'$
whose vertices $v_0, v_1, \ldots, v_{2\ell+2}$ are
\[
v_i=
\begin{cases}
(2,1,0,0,\ldots,0), \qquad \qquad &i=1, \\
(0,2,1,0,\ldots,0), &i=2, \\
(0,\ldots,0,\underbrace{1}_{i\text{-th}},\underbrace{1}_{{(i+1)}\text{-th}},0,\ldots,0) &i=3,\ldots,2l+1, \\
(1,0,\ldots,0,1), &i=2l+2, \\
(0,\ldots,0), &i=0.
\end{cases}
\]
As usual, a routine computation says that $\vol(\Pc)=3/d'!$. 
Let $v \in \RR^{d'+1}$ be the point
\begin{eqnarray*}
\frac{1}{3}\left\{(v_0,1)+(v_1,1)+(v_2,1)\right\}
+\frac{1}{3}\sum_{q=2}^{\ell+1}(v_{2q},1)
+\frac{2}{3}\sum_{q=2}^{\ell+1}(v_{2q-1},1)
\end{eqnarray*}
belonging to $\RR^{d'}$.  Then
\[
v = (1,1,\ldots,1,\ell+1) \in \ZZ^{d'+1}.
\]
Thus Lemma \ref{computation} guarantees that
$\d_{\ell+1} \geq 1$ and $\d_{\ell+1}^* \geq 1$.
Hence $\d_{\ell+1}=\d_{\ell+2}=1$,
as required.

\medskip

\noindent
{\bf (Second Step)}
Let $k \geq 1$.
We write $\Pc \subset \RR^{d'}$ for the integral simplex 
of dimension $d'$
with the vertices $v_0,v_1,\ldots,v_{3k+2\ell+2}$ as follows:
\begin{itemize}
\item
$v_0 = (0, 0, \ldots, 0)$,
\item
$v_1 = (1,1,1,0,0,\ldots,\underbrace{0}_{(3k+2)\text{-th}},1,1,\ldots,1)$,
\item
$v_2 = (0,1,1,1,0,\ldots,\underbrace{0}_{(3k+2)\text{-th}},1,1,\ldots,1)$,
\item
$v_i =
(0,\ldots,0,\underbrace{1}_{i\text{-th}},\underbrace{1}_{(i+1)\text{-th}},\underbrace{1}_{(i+2)\text{-th}},0,0,\ldots,
\underbrace{0}_{(3k+2)\text{-th}},1,0,1,0,\ldots,1,0)$
\newline 
\hfill for $i=3,4,5,\ldots,3k$, 
\item
$v_{3k+1} = (1,0,0,\ldots,0,1,\underbrace{1}_{(3k+2)\text{-th}},1,0,1,0,\ldots,1,0)$,
\item
$v_{3k+2} = (1,1,0,0,\ldots,0,\underbrace{1}_{(3k+2)\text{-th}},1,0,1,0,\ldots,1,0)$,
\item
$v_i = (0,0,\ldots,\underbrace{0}_{(3k+2)\text{-th}},
\ldots,0,\underbrace{1}_{i\text{-th}},0,1,0,\ldots,1,0)$,\\
\newline 
\hfill for $i=3k+3,3k+5,\ldots,3k+2\ell+1$,
\item
$v_i = (0,0,\ldots,\underbrace{0}_{(3k+2)\text{-th}},
\ldots,0,\underbrace{1}_{i\text{-th}},1,0,1,0,\ldots,1,0)$,\\
\newline 
\hfill for $i=3k+4,3k+6,\ldots,3k+2\ell+2$.
\end{itemize}

\bigskip

\noindent
Let $A$ denote the $(3k+2) \times (3k+2)$ matrix
$$
\left|\text{\Huge{A}}\right|=
\underbrace{
\begin{vmatrix}
1               &1              &1              &0      &\cdots    &\cdots         &\cdots         &0  \\
0               &1              &1              &1      &\ddots    &               &\text{\huge{0}}&\vdots   \\
\vdots          &\ddots         &\ddots         &\ddots &\ddots    &\ddots         &               &\vdots   \\
\vdots          &               &\ddots         &\ddots &\ddots    &\ddots         &\ddots         &\vdots   \\
\vdots          &               &               &\ddots &\ddots    &\ddots         &\ddots         &0        \\
0               &               &\text{\huge{0}}&       &\ddots    &1              &1              &1        \\
1               &               &               &       &          &0              &1              &1        \\
1               &1              &0              &\cdots &\cdots    &\cdots         &0              &1
\end{vmatrix}}_{(3k+2)\times(3k+2)}.
$$

\medskip

\noindent
Then a simple computation on determinants enables us to show that
\begin{eqnarray*}
d! \vol(\Pc)=
\underbrace{
\begin{vmatrix}
& &                & & &               & \\
& &\text{\Huge{A}} & & &\text{\Huge{*}}& \\
& &                & & &               & \\
& &                & &1&               & \\
& &\text{\Huge{0}} & & &\ddots         & \\
& &                & & &               &1
\end{vmatrix}}_{(3k+2+2\ell)\times(3k+2+2\ell)}
=\left|\text{\Huge{A}}\right|=3.
\end{eqnarray*}
\newpage
\noindent
One has
$$
\frac{1}{3}\left\{ (v_0,1)+(v_1,1)+\cdots+(v_{3k+4},1)\right\}
+\frac{2}{3}\left\{ (v_{3k+5},1)+(v_{3k+7},1)+\cdots+(v_{3k+2\ell+1},1)\right\} \\
$$
$$
+\frac{1}{3}\left\{ (v_{3k+6},1)+(v_{3k+8},1)+\cdots+(v_{3k+2\ell+2},1)\right\} \\
$$
$$
=(1,\ldots,1,k+1,1,k+2,1,\ldots,k+\ell,1,k+\ell+1) \in \ZZ^{d'+1}.
$$
Hence $\d_{k+\ell+1}=\d_{2k+\ell+2}=1$, as required.
\hspace{7cm}
\end{proof}

In order to complete a proof of 
the ``If'' part
of Theorem \ref{main} when $\sum_{i=0}^{d} \delta_i = 3$,
we must show the existence of
an integral convex polytope $\Pc \subset \RR^d$ of dimension $d$ 
whose $\delta$-vector coincides with
$(1, 0, \ldots, 0, \underbrace{1}_{m\text{-th}}, 0, \ldots, 0,\underbrace{1}_{n\text{-th}}, 0, \ldots ,0)$, where
$1 < m < n < d$ and $n-m-1 \leq m-1 \leq d-n$.

First, Lemma \ref{first} says that
there exists an integral convex polytope 
whose $\delta$-vector coincides with
\[
(1, 0, \ldots, 0, \underbrace{1}_{(n-m)\text{-th}}, 0, \ldots, 0, \underbrace{1}_{(2n-2m)\text{-th}}, 0, \ldots ,0)
\in \ZZ^{3n-3m+3}.
\]

Second, Lemma \ref{second} guarantees that
there exists an integral convex polytope whose $\delta$-vector coincides with
\[(1, 0, \ldots, 0, \underbrace{1}_{m\text{-th}}, 0, \ldots, 0, \underbrace{1}_{n\text{-th}}, 0, \ldots ,0)
\in \ZZ^{n+m+3}.\]

Finally, by using Lemma 1.1, 
there exists an integral convex polytope $\Pc$ of dimension $d$ with
\[
\delta(\Pc)
= (1, 0, \ldots, 0, \underbrace{1}_{m\text{-th}}, 0, \ldots, 0, \underbrace{1}_{n\text{-th}}, 0, \ldots ,0)
\in \ZZ^{d+1},
\]
as desired.

\newpage


\bigskip

\begin{thebibliography}{9}
\bibitem{Ehrhart}
E. Ehrhart, ``Polyn\^{o}mes Arithm\'{e}tiques et
M\'{e}thode des Poly\`{e}dres en Combinatoire,''
Birkh\"{a}user, Boston/Basel/Stuttgart, 1977.

\bibitem{HibiRedBook}
T. Hibi, ``Algebraic Combinatorics on Convex Polytopes,''
Carslaw Publications, Glebe NSW, Australia, 1992.

\bibitem{HibiLBT}
T. Hibi, A lower bound theorem for Ehrhart polynomials
of convex polytopes, {\em Adv. in Math.} {\bf 105} (1994),
162 -- 165.

\bibitem{Scott}
P. R. Scott, On convex lattice polygons,
{\em Bull. Austral. Math. Soc.} {\bf 15} (1976), 395 -- 399.

\bibitem{StanleyEC}
R. P. Stanley, ``Enumerative Combinatorics, Volume 1,''
Wadsworth \& Brooks/Cole, Monterey, Calif., 1986.

\bibitem{StanleyDRCP}
R. P. Stanley, Decompositions of rational convex polytopes,
{\em Annals of Discrete Math.} {\bf 6} (1980), 333 -- 342.

\bibitem{StanleyJPAA}
R. P. Stanley,  On the Hilbert function of a graded
Cohen--Macaulay domain, {\em J. Pure and Appl. Algebra}
{\bf 73} (1991), 307 -- 314.

\end{thebibliography}
\end{document}